\documentclass[11pt]{amsart}
\usepackage{amssymb}

\usepackage[T1]{fontenc} 
\usepackage[utf8]{inputenc} 

\usepackage[unicode,colorlinks,
	citecolor=blue,
	urlcolor=blue,
	linkcolor=blue
]{hyperref}

\usepackage{soul} 

\usepackage{url}

\usepackage{mathtools}
\usepackage{amscd}
\usepackage{tikz-cd}

\theoremstyle{plain}
\newtheorem{thm}{Theorem}[section]
\newtheorem{prop}[thm]{Proposition}
\newtheorem{lem}[thm]{Lemma}
\newtheorem{cor}[thm]{Corollary}
\newtheorem{que}[thm]{Question}

\theoremstyle{definition}
\newtheorem{defn}[thm]{Definition}

\theoremstyle{remark}
\newtheorem{remark}[thm]{Remark}

\numberwithin{equation}{section}

\renewcommand{\epsilon}{\varepsilon}
\renewcommand{\phi}{\varphi}
\renewcommand{\setminus}{\smallsetminus}
\renewcommand{\preceq}{\preccurlyeq}

\newcommand{\N}{\mathbb{N}}
\newcommand{\Q}{\mathbb{Q}}

\newcommand{\Z}{\mathbb{Z}}

\newcommand{\lb}{\lbrace}
\newcommand{\lk}{\lbrack}
\newcommand{\rb}{\rbrace}
\newcommand{\rk}{\rbrack}

\newcommand{\ra}{\rightarrow}

\begin{document}

\title{Countable Ordered Groups and Weihrauch Reducibility}

\author{Ang Li}

\date{\today}

\begin{abstract}
This paper continues to study the connection between reverse mathematics and Weihrauch reducibility. In particular, we study the problems formed from Maltsev's theorem \cite{1949} on the order types of countable ordered groups. Solomon \cite{CountableOrderedGroup} showed that the theorem is equivalent to $\mathsf{\Pi_1^1}$-$\mathsf{CA_0}$, the strongest of the big five subsystems of second order arithmetic. We show that the strength of the theorem comes from having a dense linear order without endpoints in its order type. Then, we show that for the related Weihrauch problem to be strong enough to be equivalent to $\mathsf{\widehat{WF}}$ (the analog problem of $\mathsf{\Pi_1^1}$-$\mathsf{CA_0}$), an order-preserving function is necessary in the output. Without the order-preserving function, the problems are very much to the side compared to analog problems of the big five.
\end{abstract}

\maketitle

\section{Introduction}

We follow the definitions and notations used by Solomon in \cite{CountableOrderedGroup}.

\begin{defn}
    An $\emph{ordered group}$ is a pair $(G, \leq_G)$ where $G$ is a group with identity $e$, $\leq_G$ is a linear order on $G$, and for any $a, b, g\in G$, if $a\leq_G b$, then $ag \leq_G bg$ and $ga \leq_G gb$.	
\end{defn}
\par
We will suppress the subscript of the linear order when it is clear from the context.

    \begin{defn}
        Given an ordered group $G$ and a linear order $X$, $X$ is the $\emph{order type}$ of $G$ if there is an order-preserving bijection $f$ from $G$ to $X$, which is denoted by $G\cong X$.
    \end{defn}
\par
Groups $\Z$ and $\Q$ with addition are two examples of ordered groups. From now on, we will use $\Z$ and $\Q$ for the order types of these two ordered groups.
\par
In this work, we shall focus on the class of countable ordered groups. Maltsev classified the possible order types. Before mentioning his result, we need the definitions of products of linear orders and $\Z^X$. 

    \begin{defn}
        Given two linear orders $(X,\leq_X)$ and $(Y,\leq_Y)$, the \emph{product} $XY$ is the linear order $(Z,\leq_Z)$ such that $Z=\lb \langle x,y \rangle:x\in X,y\in Y\rb$ and $\langle x_0,y_0\rangle<_Z \langle x_1, y_1\rangle$ if and only if $y_0<_Y y_1$, or $y_0=y_1$ and $x_0<_X x_1$. 
    \end{defn}

    \begin{defn}
        Given a linear order $X$, $\Z^X$ is the set of functions $f$ from $X$ to $\Z$ with finite support. If $f\not=g$, then $f<g$ if and only if $f(x)<_{\Z}g(x)$ where $x$ is the maximum value of $X$ on which $f$ and $g$ disagree.
    \end{defn}
\par
For any linear order $X=\lb x_0,x_1,\ldots\rb$, we can view elements of $\Z^X$ as finite sums $\lb \sum_{i\in I}r_ix_i:r_i\in\Z\setminus\lb 0\rb,|I|\in\N\rb$. So, we can put a natural abelian group structure on $\Z^{X}$ by addition of functions. We call this the standard group structure of $\Z^X$. In general, there are many other groups of order type $\Z^X$.
\par
Maltsev \cite{1949} proved that the order type of a countable ordered group is either $\Z^{\alpha}$ or $\Z^{\alpha}\Q$ where $\alpha$ is an ordinal. For an element $\langle z_0,z_1,\ldots,q\rangle$ of $\Z^{\alpha}\Q$, we call $q$ the $\Q$-\emph{coordinate} and $z_{\beta}$ the $\beta$-\emph{coordinate} when $\beta<\alpha$. 
\par
Reverse mathematics studies subsystems of $\mathsf{Z}_2$, the system of second order arithmetic. It was initially built on a profound observation that over a relatively weak base theory $\mathsf{RCA_0}$, most theorems are equivalent to a small number of subsystems, i.e., the $\emph{big five}$: $\mathsf{RCA_0}$, $\mathsf{WKL_0}$, $\mathsf{ACA_0}$, $\mathsf{ATR_0}$, and $\mathsf{\Pi_1^1}$-$\mathsf{CA_0}$. See \cite{Simpson} for references about reverse mathematics.
\par
Solomon \cite{CountableOrderedGroup} showed that Maltsev's theorem is on the level of $\mathsf{\Pi_1^1}$-$\mathsf{CA_0}$ in reverse mathematics:
    \begin{thm}\label{type}
        The following are equivalent over $\mathsf{RCA_0}$:
        \begin{enumerate}
            \item $\mathsf{\Pi_1^1}$-$\mathsf{CA_0}$;
	   \item Let $G$ be a countable ordered group. There is a well order $\alpha$ and $\epsilon \in \lb 0,1 \rb$ such that $\Z^{\alpha} \Q^{\epsilon}$ is the order type of $G$;
	   \item Let $G$ be a countable ordered abelian group. There is a well order $\alpha$ and $\epsilon \in \lb0,1\rb$ such that $\Z^{\alpha}\Q^{\epsilon}$ is the order type of $G$.
        \end{enumerate}
    \end{thm}
\par
Statements like the ones in this theorem and many others in ``classical'' mathematics can be written as follows: $$(\forall x\in X)(\exists y\in Y)\lk\phi(x)\ra\psi(x,y)\rk.$$ We can naturally think this as a computational problem, i.e., given an input $x$ such that $\phi(x)$, we want to produce an output $y$ such that $\psi(x,y)$. Such computational problems can be represented by partial multi-valued functions $f:\subseteq X\rightrightarrows Y$, which are just relations $f\subseteq X\times Y$, such that $f(x)=\lb y\in Y:\psi(x,y)\rb$ for each $x$ that $\phi(x)$ holds. This allows us to classify the uniform computational contents of theorems using Weihrauch reducibility.
\par
To introduce Weihrauch reducibility, we need a notion of computability on Baire space. Usually, there are two ways. One of them uses Type-2 Turing machines. See Weihrauch's book on Computable Analysis \cite{Weihrauch} for references. Here, we present the other one.

    \begin{defn}
        A single-valued function $f:\subseteq \N^{\N} \ra \N^{\N}$ is \emph{computable} if there is a total computable function $g: \N^{<\N} \ra \N^{<\N}$ such that:
        \begin{itemize}
            \item $g(\sigma)\preceq g(\tau)$ when $\sigma\preceq \tau$,
            \item $f(x)=y$ if and only if for any $n$, there exists $m$ such that $y\upharpoonright n\preceq g(x\upharpoonright m)$. 
        \end{itemize}
    \end{defn}
\par
Using realizers, we can define computability for multi-valued functions.

    \begin{defn}
        A single-valued function $f:\subseteq\N^{\N} \ra \N^{\N}$ is a \emph{realizer} for a multi-valued function $g:\subseteq\N^{\N} \rightrightarrows \N^{\N}$ if $$(\forall p\in\mathrm{dom}(g))\lk f(p)\in g(p)\rk.$$
        We say that $g$ is \emph{computable} if it has a computable realizer.
    \end{defn}
\par
In order to compare the uniform computational content of problems on spaces other than Baire space, we introduce represented spaces.

    \begin{defn}
        A \emph{represented space} is a pair $(X, \delta_X)$ where $\delta_X$ is a surjection: $\subseteq\N^{\N}\ra X$. If $\delta_X(p) = x$, then we call $p$ a \emph{name} for $x$.
        \par
        A single-valued function $F$ on Baire space is a \emph{realizer} of a multi-valued function $f:\subseteq X\rightrightarrows Y$, where $X,Y$ are represented spaces, if $$(\forall p\in\mathrm{dom}(f\circ \delta_X ))\lk \delta_Y\circ F(p)\in f\circ\delta_X(p)\rk.$$ 
    \end{defn}
\par
Now, we can define Weihrauch reducibility for problems on arbitrary spaces.

    \begin{defn}
        Let $f$, $g$ be partial multi-valued functions from $X$ to $Y$ and from $Z$ to $W$. Then, $f$ is \emph{Weihrauch reducible} to $g$, denoted $f\leq_{\sf{W}} g$ if there are computable $\Phi$, $\Psi$ on Baire space such that for any realizer $G$ of $g$, $\Psi\circ\langle\text{Id},G\circ \Phi\rangle$ is a realizer of $f$. Equivalently, if $p$ is a name of some $x\in\mathrm{dom}(f)$:
        \begin{itemize}
            \item $\Phi(p)$ is a name of some $z\in\mathrm{dom}(g)$, and
            \item given a name $q=G\circ\Phi(p)$ for some element of $g(z)$, $\Psi(p,q)$ is a name for some element in $f(x)$.
        \end{itemize}
        $\Phi,\Psi$ are called \emph{forward functional} and \emph{backward functional}, respectively.
    \end{defn}
\par    
Weihrauch reducibility is reflexive and transitive. Therefore, we can define the equivalence relation $\equiv_{\mathsf{W}}$ by $f\equiv_{\mathsf{W}}g$ if $f\leq_{\mathsf{W}}g$ and $g\leq_{\mathsf{W}}f$. The equivalence classes of $\equiv_{\sf{W}}$ are called \emph{Weihrauch degrees}.
\par
Mathematical problems can be combined in many natural ways to form new problems. Here, we give a very short list of some algebraic operations on problems that are relevant to this paper. By $X^{\ast}:=\bigcup^{\infty}
_{i=0}(\lb i\rb\times X^i)$ we denote the set of words over $X$, where $X^i:=\bigtimes^i_{j=1} X$ is the $i$-fold Cartesian product and $X^0$ is a singleton set with the empty tuple as the only element. See \cite{survey} for more about Weihrauch degree, represented spaces, and algebraic operations.

\begin{defn}
    Let $f :\subseteq X\rightrightarrows Y$, $g :\subseteq Z \rightrightarrows W$, and $h :\subseteq Y \rightrightarrows Z$ be multi-valued functions. We define the following operations:
    \begin{enumerate}
        \item composition $h\circ f:\subseteq X \rightrightarrows Z$, $(h\circ f)(x):=\lb z\in Z: (\exists y\in f(x))\lk z\in h(y)\rk\rb$ when $x\in\mathrm{dom}(h\circ f):=\lb x\in\mathrm{dom}(f):f(x)\subseteq \mathrm{dom}(h)\rb$;
        \item product $f\times g:\subseteq X\times Z\rightrightarrows Y\times W$, 
            
        $(f\times g)(x,z):= f(x)\times g(z)$, and 
            
        $\mathrm{dom}(f\times g):= \mathrm{dom}(f)\times\mathrm{dom}(g)$; 
        \item finite parallelization $f^{\ast}:\subseteq X^{\ast}\rightrightarrows Y^{\ast}$, 
            
        $f^{\ast}(i, x) := \lb i\rb\times f^i(x)$, and 
            
        $\mathrm{dom}(f^{\ast}):=\mathrm{dom}(f)^{\ast}$;
        \item parallelization $\widehat{f}:\subseteq X^{\N}\rightrightarrows Y^{\N}$, 
            
        $\widehat{f}(\bigtimes_i x_i) := \bigtimes_{i\in\N} f(x_i)$, and 
            
        $\mathrm{dom}(\widehat{f}):=\mathrm{dom}(f)^{\N}$.
    \end{enumerate}
\end{defn}
\par
Here, (2), (3), and (4) are operations that can be lifted to the Weihrauch degrees. For (1), Brattka and Pauly \cite{algebraic} showed that $f\ast g:=\max_{\leq_{\mathsf{W}}}\lb f_0\circ g_0:f_0\leq_{\mathsf{W}}f,g_0\leq_{\mathsf{W}}g\rb$ exists for any problems $f,g$, which is called the \emph{compositional product}.
\par
Much work has been done previously on connecting reverse mathematics and Weihrauch degree, which was started by Gherardi and Marcone \cite{bridge}. The big five systems have been identified with some problems in the Weihrauch lattice as their counterparts, e.g.\ $\mathsf{RCA_0}$ to computable problems, $\mathsf{WKL_0}$ to the problem $\mathsf{C_{2^{\N}}}$, and $\mathsf{ACA_0}$ to iterations of the problem $\mathsf{lim}$. However, it is more complicated to come up with an analog for $\mathsf{ATR_0}$. One of the problems identified is $\mathsf{C_{\N^{\N}}}$. Meanwhile, there is a natural Weihrauch degree $\mathsf{\widehat{WF}}$ for $\mathsf{\Pi_1^1}$-$\mathsf{CA_0}$. Here, $\mathsf{WF}$ is the problem that given a tree $T\subseteq\N^{<\N}$, the output is $0$ if $T$ is well-founded and $1$ otherwise. Cipriani, Marcone, Valenti \cite{manlio} first investigated this connection for a theorem equivalent to the system $\mathsf{\Pi_1^1}$-$\mathsf{CA_0}$. We continue this study for countable ordered groups.
\par
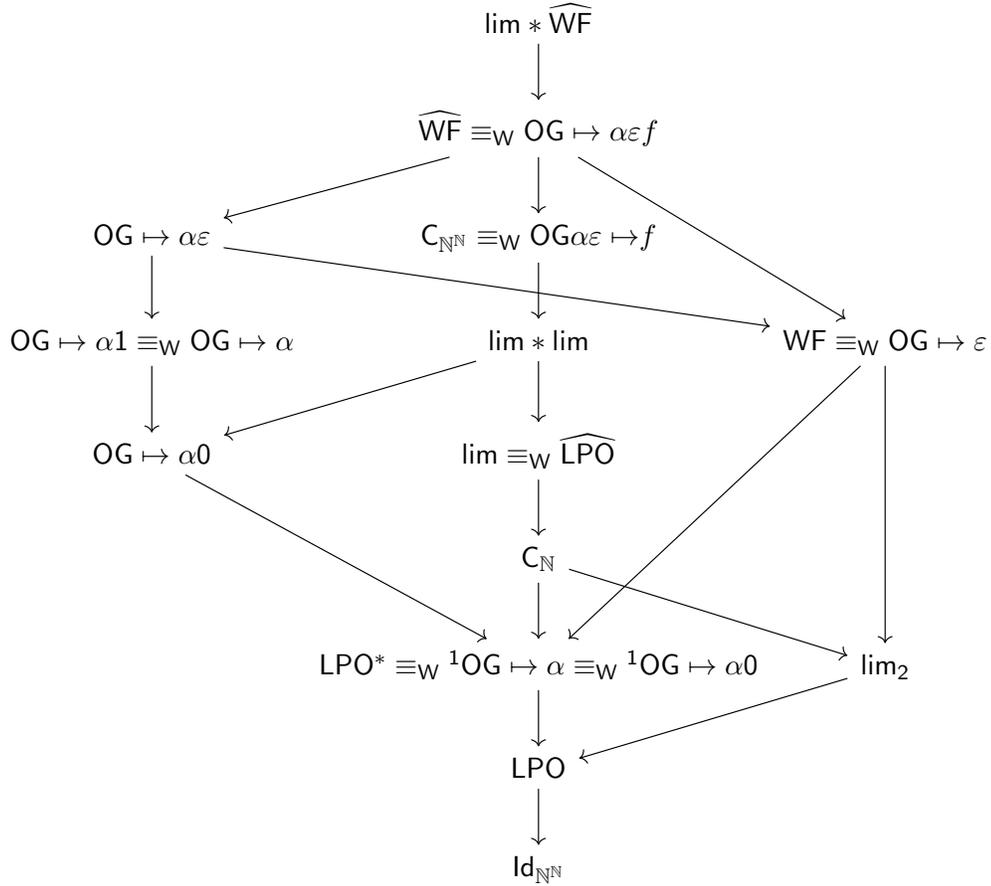
\begin{figure}[h!]
    \centering
    \begin{tikzcd}[column sep=0in, row sep=0.3in]
        & & \mathsf{\widehat{WF}}\equiv_{\sf{W}}\mathsf{OG\mapsto\alpha\epsilon}f \arrow[dl] \arrow[d] \arrow[ddr] & & \\
        & \mathsf{OG\mapsto\alpha\epsilon} \arrow[d] \arrow[drr] & \mathsf{C_{\N^{\N}}}\arrow[d] & & \\
        & \mathsf{OG\mapsto\alpha1}\equiv_{\sf{W}}\mathsf{OG\mapsto\alpha} \arrow[d] &  \mathsf{lim}\equiv_{\sf{W}}\mathsf{\widehat{LPO}} \arrow[d] \arrow[dl] & \mathsf{WF}\equiv_{\sf{W}}\mathsf{OG\mapsto\epsilon}  \arrow[dd] & \\
        & \mathsf{OG\mapsto\alpha0} \arrow[dr] & \mathsf{C}_{\N} \arrow[d] \arrow[dr] & & \\
        & & \mathsf{LPO^{\ast}}\equiv_{\sf{W}}\mathsf{^1OG\mapsto\alpha}\equiv_{\sf{W}}\mathsf{^1OG\mapsto\alpha0} \arrow[d] & \mathsf{lim_2} \arrow[dl] & \\
        & & \mathsf{LPO} \arrow[d] & & \\
        & & \mathsf{Id_{\N^{\N}}} & &
    \end{tikzcd}
    \caption{The arrows in this diagram are strict. When there are no arrows between two degrees, it means that those two degrees are incomparable.}
    \label{fig:diagram1}
\end{figure}

Figure \ref{fig:diagram1} summarizes the results in this work. Definitions of Weihrauch degrees in the diagram will be given later.

\section{Countable Ordered Groups}

Usually, for one statement, there are multiple ways to frame it as a problem. For the statements in Theorem \ref{type}, we could either let the input be any countable ordered group or restrict it to the abelian case. We have various choices for what information we output. 
\par
First, we need to clarify how we code the information of the output. There is no canonical way to code an ordinal $\alpha$. So, we code $\alpha$ by a set $A\subseteq\N$ and a relation $<_A\subseteq A\times A$ such that $<_A$ well-orders $A$ and this well-order is isomorphic to $\alpha$. We call $(A,\leq_A)$ an $\omega$-copy of $\alpha$. Given an $\omega$-copy of $\alpha$, we can find an $\omega$-presentation of $\Z^{\alpha}\Q^{\epsilon}$. We can also code the order-preserving functions in Baire Space. Now, we make a list of the problems.
\begin{itemize}
    \item $\mathsf{OG\mapsto\alpha\epsilon}$: given a countable ordered group $\langle G,\leq\rangle$, output a copy of the well-order $\alpha$ and $\epsilon=0,1$ such that the order type of $G$ is $\Z^{\alpha}\Q^{\epsilon}$.
    \item $\mathsf{OG\mapsto\alpha\epsilon}f$: the output also includes a presentation of an order-preserving function $f$ from $G$ to $\Z^{\alpha}\Q^{\epsilon}$.
    \item $\mathsf{AOG\mapsto\alpha\epsilon}f$: the input group has to be abelian.
    \item $\mathsf{OG\mapsto\alpha}$: only output a copy of the well-order $\alpha$.
    \item $\mathsf{OG\mapsto\epsilon}$: only output the one bit $\epsilon$.
\end{itemize}
\par 
We did not include more problems that require the input to be abelian because it does not matter. This is not surprising, as it is also the case in reverse mathematics. 
\par
In order to understand these problems, we need to understand the complexity of the well-order $\alpha$ in the order types of all countable ordered groups. First, we introduce a few notions from ordered group theory that are useful to us. For more about ordered groups, see Kokorin and Kopytov \cite{kokorin}.
\par
Notice that given a quotient group $G/H$, it is not guaranteed that it will inherit an order when $G$ is ordered. The following class of normal subgroups ensures the orderability of the quotient group. 
\begin{defn}
    A normal subgroup $H$ of an ordered group $G$ is \emph{convex} if for all $a,b\in H$ and $g\in G$, if $a\leq g\leq b$ then $g\in H$. The induced order on $G/H$ is defined as follows: $$aH\leq_{G/H} bH\leftrightarrow (aH=bH)\vee (aH\not =bH\wedge a<_G b).$$
\end{defn}
\begin{defn}
    For an ordered group G with identity $e$, $|x|=x$ if $x>_Ge$ (or equivalently, $x$ is a \emph{positive} element) and $|x|= x^{-1}$ otherwise. For any $a,b\in G$, $a$ is \emph{Archimedean less} than $b$, denoted $a\ll b$, if $|a^n|<|b|$ for any $n\in\N$. If there exist $n,m\in\N$ such that $|a^n|\geq|b|$ and $|b^m|\geq|a|$, then $a$ and $b$ are \emph{Archimedean equivalent}, denoted $a\approx b$.   
\end{defn}
\begin{defn}
    Let $G$ be an ordered group. The set $\mathrm{Arch(G)}$ is a set of unique representatives of the Archimedean classes of $G$. $$\mathrm{Arch}(G)=\lb g \in G :(\forall h\in G)\lk(h<_{\N} g)\ra \neg(h\approx g)\rk\rb.$$ It is ordered by taking $x < y$ if and only if $x\ll y$. We also define $W(\mathrm{Arch}(G))$ to be the largest initial segment of $\mathrm{Arch}(G)$ that is well-ordered.
\end{defn}
\begin{lem}\label{Arch}
    If the order type of $G$ is $\Z^{\alpha}$, $\mathrm{Arch(G)}$ is a copy of $\alpha$. If the order type of $G$ is $\Z^{\alpha}\Q$, $W(\mathrm{Arch}(G))$ is a copy of an ordinal $\beta\geq\alpha$.
\end{lem}

\begin{proof}
    By the proof of Theorem 1 in Chapter 7 in \cite{kokorin}, there exists a unique chain of convex normal subgroups $\lb e\rb=A_0\subset A_1\subset \ldots\subset A_{\gamma}\subset A_{\gamma+1}\subset\ldots\subset A_{\alpha}\subseteq G$ that has the following properties: $A_{\gamma+1}/A_{\gamma}$ is an infinite cyclic group; if $\gamma$ is a limit ordinal, then $A_{\gamma}=\bigcup_{\tau<\gamma}A_{\tau}$; $G/A_{\alpha}$ has order type $\Q$ or $G=A_{\alpha}$. Suppose $b_{\gamma}'$ is the positive generator of $A_{\gamma+1}/A_{\gamma}$ for each $\gamma$. Let $b_{\gamma}$ be a positive element in $G$ such that $b_{\gamma}A_{\gamma}=b_{\gamma}'$. We claim that $(\lb b_{\gamma}:\gamma<\alpha\rb,\ll)$ is isomorphic to $\alpha$, and $\alpha$ is no more than the ordinal $\beta$ isomorphic to $W(\mathrm{Arch}(G))$. To prove the claim, we verify the following properties.
    \begin{enumerate}
        \item For any $\gamma<\tau<\alpha$, $b_{\gamma}\ll b_{\tau}$. It suffices to show that for any $\tau>0$, $b_{\tau}$ is Archimedean more than any element in $A_{\tau}$. Suppose $b_{\tau}$ is not Archimedean above a positive element $a\in A_{\tau}$. There exists $n$ such that $b_{\tau}\leq a^n$. By the convexity of $A_{\tau}$, $b_{\tau}\in A_{\tau}$. Contradiction.
        \item For any $\gamma<\alpha$ and any element $c\in A_{\gamma+1}$, $c$ is Archimedean less than or Archimedean equivalent to $b_{\gamma}$. We prove this by induction. If $c$ is in the coset $eA_{\gamma}$, then $c\ll b_{\gamma}$ since there exist $\tau<\gamma$ such that $c\ll b_{\tau}$ or $c\approx b_{\tau}$ and by the first property, $b_{\tau}\ll b_{\gamma}$ ($c=e\ll b_0$ when $\gamma=0$). If $c$ is in other cosets of $A_{\gamma+1}/A_{\gamma}$, then $c=b_{\gamma}^nd$ for some nonzero integer $n$ and some $d\in A_{\gamma}$ since $A_{\gamma+1}/A_{\gamma}$ is cyclic. Since $d$ is Archimedean less than $b_{\gamma}$ ($d=e\ll b_0$ when $\gamma=0$), we conclude that $c$ is Archimedean equivalent to $b_{\gamma}$ in this case.
        \item For any $\gamma+1<\alpha$, there does not exist an element $g\in G$ such that $b_{\gamma}\ll g\ll b_{\gamma+1}$. Assume otherwise. By the convexity of $A_{\gamma+2}$, $g\ll b_{\gamma+1}$ implies that $g\in A_{\gamma+2}$. Then, $g$ has to be in the coset $eA_{\gamma+1}$ of $A_{\gamma+2}/A_{\gamma+1}$ by the proof of the second property. So, $g\in A_{\gamma+1}$. But any element in $A_{\gamma+1}$ is not Archimedean more than $b_{\gamma}$. Contradiction.
        \item For any limit ordinal $\gamma<\alpha$ and $g\in G$ that is Archimedean less than $ b_{\gamma}$, there exists $\tau<\gamma$ such that $g\ll b_{\tau}$. Notice that $g\ll b_{\gamma}$ implies $g\in A_{\gamma+1}$ by convexity. By the proof of the second property, $g\in A_{\gamma}$. Therefore, there exists $\tau<\gamma$ such that $g\in A_{\tau}$, which implies $g\ll b_{\tau}$.
    \end{enumerate}
    By property (1), $(\lb b_{\gamma}:\gamma<\alpha\rb,\ll)$ is isomorphic to $\alpha$. Properties (3) and (4) show that $(\lb b_{\gamma}:\gamma<\alpha\rb,\ll)$ is isomorphic to an initial  segment $\mathrm{Arch}(A_{\alpha})$ of $W(\mathrm{Arch}(G))$. Therefore, $\alpha\leq\beta$. When the order type of $G$ is $\Z^{\alpha}$, $\mathrm{Arch}(G)$ is isomorphic to $\alpha$ since $G=A_{\alpha}$.
\end{proof}

\begin{lem}\label{G}
    If $G$ is an ordered group, then we can uniformly compute $\mathrm{Arch}(G)$ from $G'$.
\end{lem}

\begin{proof}
    Say $G$'s elements are listed as $\lb g_i\rb_{i\in\N}$. We first build a sequence $\lb p_i\rb_{i\in\omega}$ of infinite strings in $\N^{\N}$. We let $p_i(\langle n,m\rangle)$ be $1$ if $|g_n^i|<|g_m|$, and $0$ otherwise. Note that $p=\lim_{i\ra\infty}p_i$ exists and is computable from $G'$. Using $p$, we can compute a copy of $\mathrm{Arch(G)}$ by putting $g_n$ in when $p(\langle n,m\rangle) =1$ or $p(\langle m,n\rangle) =1$ for each $m<n$.
\end{proof}

\par
Recall that $\omega_1^{\mathrm{CK}}$ is the first non-computable ordinal. It does not have a hyperarithmetical copy \cite{spector}.

\begin{prop}\label{q}
    If a computable ordered group has order type $\Z^{\alpha}\Q$, then $\alpha\leq \omega_1^{\mathrm{CK}}$. 
\end{prop}

\begin{proof}
    By Lemma \ref{G}, we have a $0'$-computable copy of $\mathrm{Arch}(G)$. By Lemma \ref{Arch}, $\alpha$ is a well-ordered initial segment of $\mathrm{Arch}(G)$. Notice that $\omega_1^{\mathrm{CK}}+1$ cannot be an initial segment of $\mathrm{Arch}(G)$. Otherwise, $\omega_1^{\mathrm{CK}}$ is $0'$-computable. Therefore, $\alpha\leq\omega_1^{\mathrm{CK}}$. 
\end{proof}

\begin{prop}\label{main}
    If a computable ordered group has order type $\Z^{\alpha}$, then $\alpha<\omega_1^{\mathrm{CK}}$.
\end{prop}
\begin{proof}
    By Lemmas \ref{Arch} and \ref{G}, $\mathrm{Arch}(G)$ has a $0'$-computable copy, which also is a copy of $\alpha$. Therefore, $\alpha<w_1^{\mathrm{CK}}$.   
\end{proof}

\par
In reverse mathematics, the subsystem $\mathsf{ATR_0}$ does not imply Maltsev's theorem because the order type could contain a $\Q$-part. Now, we see that this complexity is also reflected in the complexity of the well-order $\alpha$ when the order type contains a $\Q$-part by showing that $\alpha$ can be $\omega_1^{\mathrm{CK}}$ in this case.

\begin{prop}\label{2way}
    There exists a computable countable ordered group with order type $\Z^{\omega_1^{\mathrm{CK}}}\Q$.
\end{prop}

\begin{proof}
    Recall that the Harrison linear order $\mathcal{H}$ is a computable linear order of order type $\omega_1^{\mathrm{CK}}(1+\Q)$. Then, the standard countable ordered group of order type $\Z^{\mathcal{H}}$ has a computable copy. Notice that $\Q$ can be embedded into $\Z^{\mathcal{H}}$. Therefore, the order type of the group is $\Z^{\alpha}\Q$ for some ordinal $\alpha$. This ordinal is at least $\omega_1^{\mathrm{CK}}$ since $\Z^{\mathcal{H}}=\Z^{\omega_1^{\mathrm{CK}}}\Z^{\omega_1^{\mathrm{CK}}\Q}$ has an interval centered at the identity element of the group that has order type $\Z^{\omega_1^{\mathrm{CK}}}$. By Proposition \ref{q}, $\alpha$ is not larger than $\omega_1^{\mathrm{CK}}$.
\end{proof}

There is another way to see this. First, we introduce the notions of trees and well-foundedness formally. Let $\N^{<\N}$ denote the set of finite strings of natural numbers, where $\lambda$ is the empty string. Given finite strings $\sigma$, $\tau$, the concatenation is denoted as $\sigma^{\smallfrown}\tau$ and the length of $\sigma$ is denoted $|\sigma|$. We write $\tau\preceq\sigma$ if $\sigma$ extends $\tau$. We also use $|\cdot|$ for sizes of sets.

\begin{defn}
    A \emph{tree} $T$ in Baire space is a nonempty subset of $\N^{<\N}$ closed under initial segments. A string $\sigma\in T$ is called a \emph{node}. We say it is on \emph{level} $l$ when $|\sigma|=l$. A function $f\in\N^{\N}$ is called a \emph{path} through $T$ if $f\upharpoonright n=f(0)^{\smallfrown}\cdots^{\smallfrown}f(n-1)\in T$ for every $n$. A tree $T$ is \emph{ill-founded} if $T$ has at least one path and \emph{well-founded} otherwise.
\end{defn}

\begin{remark}\label{gandy}
    The other way to show Proposition \ref{2way} is by contradiction using the Gandy Basis Theorem and the reduction from $\widehat{\sf{WF}}$ to $\mathsf{OG\mapsto\alpha\epsilon}f$, which we will prove later in Proposition \ref{outputeverything}. Notice that there is an effective list $\lb T_i\rb_{i\in\N}$ of computable trees in the sense that given any computable tree $T$, there is some $T_i$ such that $T$ and $T_i$ share the same infinite paths. We can make such a list by constructing $T_i$ for each c.e.\ set $W_i$ as follows: put $\sigma$ into $T_i$ if and only if there is no prefix of $\sigma$ in $W_{i,|\sigma|}$. Notice that $\text{WF}_0=\lb i:T_i\text{ is well-founded}\rb$ is $\Pi^1_1$-complete. So, there is some computable input of $\widehat{\sf{WF}}$ such that its output computes $\mathcal{O}$. Given a computable input to $\widehat{\sf{WF}}$, we get a computable countable ordered group from the forward functional of the reduction from $\widehat{\sf{WF}}$ to $\mathsf{OG\mapsto\alpha\epsilon}f$. Recall that the Gandy Basis Theorem says that if a non-empty set $A\subseteq\N^{\N}$ is $\Sigma_1^1$, then $A$ contains an element $x$ such that $\omega_1^x = \omega_1^{\mathrm{CK}}$ and $x <_T \mathcal{O}$. If $\alpha$ cannot be non-computable, then the set of order-preserving bijections from the group to its order type as a subset of Baire space is $\Pi^0_2$. So, there is an element of this set that cannot compute Kleene's $\mathcal{O}$. Then, the backward functional cannot compute $\mathcal{O}$ using $\alpha,\epsilon, f$. 
\end{remark}
\par
Notice that the proofs of the above three propositions can be relativized. Therefore, we have the following corollary.

\begin{cor}\label{relativization}
    If an $X$-computable ordered group has order type $\Z^{\alpha}$, then $\alpha$ is $X$-computable. If an $X$-computable ordered group has order type $\Z^{\alpha}\Q$, then $\alpha\leq\omega_1^X$. There exists an $X$-computable countable ordered group with order type $\Z^{\omega_1^{X}}\Q$.
\end{cor}
\par
Next, we define the Kleene-Brouwer order. Then, we make use of it to show how difficult it is to decide the one-bit information $\epsilon$. 

\begin{defn}
    The \emph{Kleene--Brouwer order} $\leq_{\mathrm{KB}}$ on $\N^{<\N}$ is as follows: given any $\sigma$, $\tau\in\N^{<\N}$, $\sigma\leq_{\mathrm{KB}}\tau$ if and only if $\tau\preceq\sigma$ or there is some $j<\min\lb|\sigma|,|\tau|\rb$ with $\sigma(j)<\tau(j)$ and $\sigma(i)=\tau(i)$ for all $i<j$. 
\end{defn}

\par
Given a tree $T$, let $\mathrm{KB}(T)$ be the restriction of the Kleene--Brouwer order to $T$. It is a linear order that is a well-order if and only if $T$ is well-founded.

\begin{prop}\label{onepoint}
    $\sf{OG\mapsto\epsilon} \equiv_{\sf{W}}\sf{WF}$.
\end{prop}

\begin{proof}
    ``$\leq$'' Given a countable ordered group $G$, we build a tree $T$ by trying to embed the rationals into the order. We fix an enumeration of all rational numbers $\lb q_i\rb_{i\in\N}$. We define $T$ as follows: any $\sigma$ is in $T$ if and only if the map from $q_i$ to $\sigma(i)\in G$ for $i<|\sigma|$ preserves the order. Then, $T$ is well-founded if and only if $\epsilon=0$.
    \par 
    ``$\geq$'' Given a tree $T$, let $\mathrm{KB}(T)=\lb x_0,x_1,\ldots\rb$. We can consider the standard group structure of order type $\Z^{\mathrm{KB}(T)}$. By \cite[Lemma 4.5]{CountableOrderedGroup}, $\mathrm{KB}(T)$ is a well-order if and only if $\epsilon=0$ for $G$. So, $T$ is well-founded if and only if $\epsilon=0$ for $G$.
\end{proof}
\par
We are prepared to show the following proposition.

\begin{prop}\label{outputeverything}
    $\mathsf{OG\mapsto\alpha\epsilon}f\equiv_{\sf{W}}\widehat{\sf{WF}}$.
\end{prop}

\begin{proof}
    \par
    ``$\leq$'' We prove that $\widehat{\sf{WF}}\geq_{\sf{W}} \mathsf{OG\mapsto\alpha\epsilon}f$. First, we fix a list of rational numbers $\lb q_l\rb_{l\in\N}$. We may assume that the input $X$ is an $\omega$-copy of $\langle G,\leq\rangle$ where $G$ has order type $\Z^{\alpha}\Q^{\epsilon}$. Similar to Remark \ref{gandy}, the forward functional can take in an $X$-computable input and output a sequence of trees such that the output of $\mathsf{\widehat{WF}}$ is $\Pi^1_1$-complete relative to $X$.
    \par
    As in Proposition \ref{onepoint}, one application of $\sf{WF}$ is enough to output $\epsilon$. We assume that $\epsilon=1$ since for the case $\epsilon=0$, the proof below can be adapted to output $\alpha$ and $f$. Now, we construct the backward functional $\Psi$ so that it outputs $\alpha$ and an order-preserving bijection from $G$ to $\Z^{\alpha}\Q$ using answers to some $\Pi^1_1$ questions and the original input. Given the ordered the group $\langle G,\leq\rangle$ with enumeration $\lb g_i\rb_{i\in\N}$, it can locate the identity element $e\in G$. So, let us assume that $g_0=e$. Then, we show that $\Psi$ can tell which copy of $\Z^{\alpha}$ a positive element is in and guess the ordinal $\alpha$ with the ordering-preserving map. Suppose this has been done for elements $g_0, g_1,\ldots,g_{n-1}$. Also, $\Psi$ can tell whether $g_n>e$ or not. We can skip the negative elements since their inverses will appear and determine where they should be mapped to. So, we can assume that $g_n$ is positive. For each pair $(g_i,g_n)$ for $i<n$, $\Psi$ can use the answer to the question whether $\eta$ can be embedded in between $g_i,g_n$. Answer no means that $g_i,g_n$ should be mapped to the same $\Z^{\alpha}$ copy (equivalently, share the same $\Q$-coordinate in $\Z^{\alpha}\Q$). When the answers are all yes, the backward functional check whether $g_i<g_n<g_j$ for $i,j<n$ or $g_n>g_i$ for all $i<n$. Say $g_i$'s $\Q$-coordinate is $q_{l_i}$ for $i<n$. The backward functional chooses the  $q_l$ between $q_{l_i},q_{l_j}$ such that $l$ is the least in the first case, and $q_l>q_{l_i}$ with the least $l$ for all $i<n$ in the second case to be the $\Q$-coordinate of $g_n$. Also, all the other coordinates of $g_n$ will be $0$. The choice of $q_l$ is necessary to ensure the bijectivity of the order-preserving map that is going be produced. If one of the trees is well-founded, say the one for the pair $(g_{i_0},g_n)$, the backward function assigns $g_n$ $\Q$-coordinate $q_{l_{i_0}}$. 
    \par
    Suppose we already have the current guess $\sum_{r=0}^{m-1} \alpha_r$ of $\alpha$, the current order-preserving function $f_{n-1}$ from each $g_i<g_n$ with $i<n$ to an element in the corresponding copy of $\Z^{\alpha}$, and $g'$ is the first element appeared in the copy of $\Z^{\alpha}$ $g_n$ is in. The backward functional can find a copy of the least ordinal $\beta$ such that the interval of elements between $g'$ and $g_n$ can be embedded into $\Z^{\beta}$. It can also compare this ordinal with $\sum_{r=0}^{m-1} \alpha_r$. If it is larger, we can find and cut off part of this copy that is isomorphic to $\sum_{r=0}^{m-1} \alpha_r$ to get a copy of a new ordinal. Let $\alpha_m$ be this ordinal and update the guess of $\alpha$ to be $\sum_{r=0}^{m} \alpha_r$. If not, the backward functional does not update its guess of $\alpha$. Notice that any $\alpha_r$ guessed is $X$-computable by Corollary \ref{relativization}. Therefore, the output of $\mathsf{\widehat{WF}}$ can be used by the backward functional to update the guesses. Next, the backward functional will update the order-preserving function to $f_n$ by mapping $g_n$ to an element in $\Z^{\sum_{r=0}^{m} \alpha_r}\Q$ such that $f_n$ is order-preserving and it can be extended to a map that is bijective from the interval between $g'$ and $g_n$ to the interval between $f_n(g')$ and $f_n(g_n)$. Eventually, the backward functional outputs $\sum_{i\in\N}\alpha_i$ and $\bigcup_{n\in\N} f_n$.
    \par
    ``$\geq$'' We prove that $\sf{AOG\mapsto\alpha\epsilon}f\geq_{\sf{W}}\widehat{\sf{WF}}$. This direction is easier since we can make use of the proof of Theorem \ref{type}. Given countably many trees $\lb T_i\rb_{i\in\N}$, define $T=\lb\lambda\rb\cup\lb i\ast T_i\rb$, where $i\ast T_i=\lb i^{\smallfrown}\sigma:\sigma\in T_i\rb$. Let $G$ be the free abelian group on the generators $a_{\tau}$, for $\tau\in T$. Order the generators by $a_{\tau}\ll a_{\gamma}$ if and only if $\tau<_{\mathrm{KB}}\gamma$. Then, by \cite[Lemma 4.5]{CountableOrderedGroup}, $\mathrm{KB}(T)$ is well-ordered if and only if $\epsilon=0$ for $G$. By the proof of \cite[Lemma 4.12]{CountableOrderedGroup}, all $T_i$'s are well-founded if and only if $\epsilon=0$, and when $\epsilon=1$, or equivalently, at least one of $T_i$ is ill-founded, $T_i$ is well-founded if and only if the $\Q$-coordinates of $f(a_{i-1})$ and $f(a_i)$ ($f(e)$ and $f(a_0)$ when $i=0$) are the same. 
\end{proof}
\par
Notice that in this proof, we showed that whether the group in the input is restricted to be abelian or not does not make a difference. It is similar for the problems we are going to discuss. Therefore, we omit any discussions about abelian groups from now on.
\par
It is natural to consider the problems $\sf{OG\mapsto\alpha\epsilon}$ and $\sf{OG\mapsto\alpha}$ after showing $\mathsf{OG\mapsto\alpha\epsilon}f$ is equivalent to $\mathsf{\widehat{WF}}$. As mentioned in Figure \ref{fig:diagram1}, they are very much to the side of the analogs of the big five and $\mathsf{C_{\N^\N}}$ is not Weihrauch above them, where $\mathsf{C_{\N^\N}}$ is Weihrauch equivalent to the problem that given an ill-founded tree in Baire space, it computes an path through it. One way to prove this is to consider the intermediate problem $\mathsf{OG\alpha\epsilon\mapsto}f$: in addition to a copy of an ordered group $\langle G,\leq\rangle$, we are given a copy of $\alpha$ and the one bit $\epsilon$ in the order type $\Z^{\alpha}\Q^{\epsilon}$ of $G$ as input, and the output is an order-preserving function $f$. 

\begin{prop}
    $\mathsf{OG\alpha\epsilon\mapsto}f\leq_{\sf{W}}\sf{C}_{\N^{\N}}$.
\end{prop}

\begin{proof}
    Given an $\omega$-copy of the group with $\epsilon$ and a copy of $\alpha$, the forward functional builds a tree $T$ in Baire space. We identify group elements and elements of $\Z^{\alpha}\Q^{\epsilon}$ with natural numbers in their corresponding $\omega$-copies. For any $\sigma\in T$ with $|\sigma|=2n$ for some $n\in\N$, $\sigma^{\smallfrown} x$ is in $T$ for an element $x$ of $\Z^{\alpha}\Q^{\epsilon}$ when the extension of the map determined by $\sigma$ that maps the $n$th element of the group to $x$ still preserves the order. When $|\sigma|$ is $2n+1$, $\sigma^{\smallfrown} g$ is in $T$ for an element $g$ of the group when the extended map that maps the $n$th element of $\Z^{\alpha}\Q^{\epsilon}$ to $g$ preserves the order. The backward functional outputs the map determined by the output path, which is a desired $f$.
\end{proof}

It is not known whether the reduction is strict.

\begin{que}
    Is $\mathsf{OG\alpha\epsilon\mapsto}f\geq_{\sf{W}}\sf{C}_{\N^{\N}}$? If not, what can we say about $\mathsf{OG\alpha\epsilon\mapsto}f$ in the Weihrauch lattice?
\end{que}

\begin{cor}\label{1}
    $\sf{OG\mapsto\alpha\epsilon} \not\leq_{\sf{W}}\sf{C}_{\N^{\N}}\ast\sf{WF}$.
\end{cor}

\begin{proof}
    Brattka, de Brecht, and Pauly \cite{Brattka} showed that $\sf{C}_{\N^{\N}}$ is closed under compositional product. Thus, $\sf{OG\mapsto\alpha\epsilon} \leq_{\sf{W}}\sf{C}_{\N^{\N}}\ast\sf{WF}$ implies $\widehat{\sf{WF}}\equiv_{\sf{W}}\mathsf{OG\mapsto\alpha\epsilon}f\leq_{\sf{W}}\mathsf{OG\alpha\epsilon\mapsto}f\ast\mathsf{OG\mapsto\alpha\epsilon} \leq_{\sf{W}}\mathsf{OG\alpha\epsilon\mapsto}f\ast\sf{C}_{\N^{\N}}\ast\sf{WF}\leq_{\sf{W}}\sf{C}_{\N^{\N}}\ast\sf{C}_{\N^{\N}}\ast\sf{WF}\equiv_{\sf{W}}\sf{C}_{\N^{\N}}\ast\sf{WF}$. So, it suffices to show $\widehat{\sf{WF}} \not\leq_{\sf{W}}\sf{C}_{\N^{\N}}\ast\sf{WF}$. 
    \par
    Assume otherwise. Let $\sf{F}$ and $\sf{G}$ be any two problems such that $\sf{F}\leq_{\sf{W}}\sf{C}_{\N^{\N}}$, $\sf{G}\leq_{\sf{W}}\sf{WF}$, and $\widehat{\sf{WF}}\leq_{\sf{W}}\sf{F}\circ\sf{G}$. Also, let $\Phi_F,\Phi_G,\Phi$ and $\Psi_F,\Psi_G,\Psi$ be forward and backward functionals for these three reductions respectively. Similar to Remark \ref{gandy}, we have some computable input $c$ for $\widehat{\sf{WF}}$ such that the output computes $\mathcal{O}$. Let $\theta\in 2$ be the output of $\sf{WF}$ given input $\Phi_G\circ\Phi(c)$. Then, the output $h=\Psi_G\langle\Phi(c),\theta\rangle$ of $\sf{G}$ is computable, which is an input for $\sf{F}$ as well. By the Gandy Basis Theorem, $\sf{C}_{\N^{\N}}$ has an output $X$ that is strictly Turing below $\mathcal{O}$ given input $\Phi_F(h)$. Then, the output $\Psi_F\langle h,X\rangle$ of $\sf{F}\circ\sf{G}$ does not compute $\mathcal{O}$, which gives a contradiction. 
\end{proof}

\begin{cor}\label{2}
    $\sf{OG\mapsto\alpha} \not\leq_{\sf{W}}\sf{C}_{\N^{\N}}$.
\end{cor}

\begin{proof}
    Notice that $\sf{OG\mapsto\alpha} \leq_{\sf{W}}\sf{C}_{\N^{\N}}$ implies $\sf{OG\mapsto\alpha\epsilon}\leq_{\sf{W}}\sf{OG\mapsto\alpha}\ast\sf{OG\mapsto\epsilon} \leq_{\sf{W}}\sf{C}_{\N^{\N}}\ast\sf{OG\mapsto\epsilon}\equiv_{\sf{W}}\sf{C}_{\N^{\N}}\ast\sf{WF}$.
\end{proof}
\par
This completes one direction of the comparison between $\sf{OG\mapsto\alpha\epsilon}$, $\sf{OG\mapsto\alpha}$ and analogs of the big five. For the other direction, we need a very weak problem called limited principle of omniscience whose finite parallelization was shown by Brattka, de Brecht, and Pauly \cite{Brattka} to be strictly below $\mathsf{C_{\N}}$: $$\mathsf{LPO}(p)=\begin{cases} 0 \text{\ if\ } (\exists n)\lk p(n)=0\rk, \\ 1 \text{\ otherwise,}\end{cases}$$ where $p\in\N^{\N}$.

\begin{prop}\label{5}
    $\sf{LPO^{\ast}}\leq_{\sf{W}}\sf{OG\mapsto\alpha}$. 
\end{prop}

\begin{proof}
    Let $(k,(p_0,\ldots,p_{k-1}))$ be any instance of $\sf{LPO^{\ast}}$. The forward functional $\Phi$ first provides $k+1$ many elements in the group $G$ so that if none of the $p_i$'s contains a zero, $G$ will be the standard group of order type $\Z^{k+1}$ and each of those $k+1$ elements will be mapped to elements of the form $(0,\ldots,0,1,0,\ldots,0)$ by an order isomorphism. If $\Phi$ finds zeros in $n$ many $p_i$'s, it will reinterpret elements in $G$ so that the order type of $G$ is $\Z^{k+1-n}$. Notice that we can first do this reinterpretation so that the order type seems to be $\Z^{k}$ instead of $\Z^{k+1}$: find a large enough finite number $l$ such that any element corresponding to $(a_0,\ldots,a_{k-1},b)$ now corresponds to $(a_0,\ldots, a_{k-2},a_{k-1}+b\cdot l)$ and the new order given by this correspondence is consistent with group multiplication. Such $l$ exists since we have only seen finitely many elements in $G$. For example, $l$ can be chosen to be $2\max|a_{k-1}|+1$, where the $\max$ is taken out of all $|a_{k-1}|$ in $(a_0,\ldots,a_{k-1},b)$ corresponding to elements that have been examined by $\Phi$. By repeating this, we can produce a group $G$ so that the order type of $G$ is $\Z^{k+1-n}$. 
    \par
    The backward functional $\Psi$ reads $(k,(p_0,\ldots,p_{k-1}))$ and a copy $c$ of $k+1-n$ at the same time. When $\Psi$ sees $k+1-n$ many elements in $c$, it  knows that there are no more than $n$ many $p_i$'s with at least one digit zero. So, $\Psi$ can read initial segments of $(k,(p_0,\ldots,p_{k-1}))$ and $c$ until the number of such $p_i$'s match exactly with what $c$ indicates. In that case, $\Psi$ knows the rest of $p_i$'s without zeros yet do not have zeros. It can output the correct answer to input $(k,(p_0,\ldots,p_{k-1}))$.  
\end{proof}
\par
Now, we want to show what $\mathsf{OG\mapsto\alpha}$ does not imply. One way to show non-reductions in the Weihrauch lattice is to find a simple witness called the first-order part of a problem. Intuitively, given a problem $\mathsf{P}$, the first-order part is a notion that captures the strongest problem with codomain $\N$ that is reducible to $\mathsf{P}$. Here, we embed $\N$ in $\N^{\N}$ by identifying $n\in\N$ with $f\in\N^{\N}$ such that $f(0)=n$ and $f(m)=0$ for all $m>0$.

\begin{defn}
    A problem $\mathsf{P}$ is \emph{first-order} if $\mathsf{P}(f)\subseteq\N$ for all $f\in dom(\mathsf{P})$. We let $\mathcal{F}$ denote the class of all first-order problems.
\end{defn}

\begin{defn}
    For any problem $\mathsf{P}$, its \emph{first-order part} $\mathsf{^1P}$ is defined as follows: $\mathsf{^1P}$ is a problem that witnesses the existence of $\max_{\leq_{\mathsf{W}}}\lb \mathsf{Q}\in\mathcal{F}:\mathsf{Q}\leq_{\mathsf{W}}\mathsf{P}\rb$.
\end{defn}
\par
This notion was first given by Dzhafarov, Solomon, and Yokoyama in \cite{firstorderpart}, where they showed that $\max_{\leq_{\mathsf{W}}}\lb \mathsf{Q}\in\mathcal{F}:\mathsf{Q}\leq_{\mathsf{W}}\mathsf{P}\rb$ exists by giving the following explicit witness.

\begin{thm}
    For a problem $\mathsf{P}$, its first-order part $\mathsf{^1P}$ can be taken to be the following problem $\mathsf{Q}$:
        \begin{itemize}
            \item the $\mathsf{Q}$-instances are all triples $\langle f, \Phi, \Psi\rangle$, where $f\in \N^{\N}$ and $\Phi$ and $\Psi$ are Turing functionals such that $\Phi(f)\in \mathrm{dom}(\mathsf{P})$ and $\Psi^{f\oplus g}(0)\downarrow$ for all $g\in \mathsf{P}(\Phi(f))$;
            \item the $\mathsf{Q}$-solutions to any such $\langle f, \Phi, \Psi\rangle$ are all $y$ such that $\Psi^{f\oplus g}(0)\downarrow= y$ for some $g\in \mathsf{P}(\Phi(f))$.
        \end{itemize}
\end{thm}
\par
See also \cite{solda} for more about first-order parts.

\begin{prop}\label{W}
    $\sf{^1OG\mapsto\alpha}\equiv_{\sf{W}}\sf{LPO^{\ast}}$.
\end{prop}

\begin{proof}
    \cite[Proposition 11.7.22]{survey} gives $\sf{Min}\equiv_{\sf{W}}\sf{LPO^{\ast}}$, where $\mathsf{Min}:\subseteq\N^{\N}\ra\N$ is the problem that outputs $\min\lb p(n):n \in\N\rb$ given input $p$. So, it suffices to show that $\sf{^1OG\mapsto\alpha}\leq_{\sf{W}}\sf{Min}$. We use the equivalent problem in the characterization of the first-order part. Given any input $\langle f,\Phi_0,\Psi_0\rangle$ of $\sf{^1OG\mapsto\alpha}$, we define the forward functional $\Phi_1$ as follows. For each triple $\langle \sigma,\tau_i,s\rangle$, $\Phi_1$ simulates $\Psi_{0,s}^{\sigma\oplus \tau_i}(0)$ where $\sigma$ is an initial segment of $f$ and $\tau_i$ is the $i$th finite string. Each time $\Psi_{0,s}^{\sigma\oplus \tau_i}(0)$ converges, $\Phi_1$ looks for the smallest finite ordinal $n$ such that it has a copy with initial segment $\tau_i$ and adds one more digit $n$ to its output. Eventually, $\Phi_1$ outputs a sequence $p$ in Baire space. Then, $\sf{Min}$ can tell us the minimum digit of this sequence.
    \par
    Let $m$ be the minimum digit. Notice that $\alpha\geq m$. Otherwise, for any copy $c$ of $\alpha$, $\Psi_0^{f\oplus c}(0)$ should converge. Then, there are initial segments $\sigma$, $\tau$ of $f$, $c$, and step $s$ such that $\Psi_{0,s}^{\sigma\oplus \tau}$ converges, which means $p$ has a digit smaller than $m$. We let the backward functional $\Psi_1$ simulate $\Psi_{0,s}^{\sigma\oplus \tau_i}(0)$ for triples $\langle \sigma,\tau_i,s\rangle$ until it converges for some $\tau_i$ such that $m$ is the smallest finite ordinal which has a copy $\gamma$ with initial segment $\tau_i$. Then, $\Psi_1$ outputs exactly $\Psi_0^{f\oplus \gamma}(0)$. Notice that $\tau_i$ can be extended to a copy $c'$ of $\alpha$. So, $\Psi_0^{f\oplus \gamma}(0)=\Psi_0^{f\oplus c'}(0)$, which is a valid output of $\sf{^1OG\mapsto\alpha}$.
\end{proof}

We only used that the output of $\sf{OG\mapsto\alpha}$ is a copy of an ordinal in this proof. So, the first-order part of any problem that outputs a copy of an ordinal is Weihrauch reducible to $\sf{LPO^{\ast}}$.
\par
Also, notice that $\sf{LPO^{\ast}}\times\sf{LPO^{\ast}}\equiv_W\sf{LPO^{\ast}}$. To see that the first order part of $\sf{^1OG\mapsto\alpha\epsilon}$ is $\sf{LPO^{\ast}}\times\sf{WF}$, we can use similar forward and backward functionals except that the problem $\sf{OG\mapsto\epsilon}$ is solved by $\sf{WF}$ and the forward functional will create two sequences for $\epsilon=0,1$ and the backward functional will use the correct minimum given $\epsilon$.
\begin{prop}\label{LW}
    $\sf{^1OG\mapsto\alpha\epsilon}\equiv_{\sf{W}}\sf{LPO^{\ast}}\times\sf{WF}$.
\end{prop}
In \cite{lim2}, Brattka, Gherardi, and Marcone showed that $\sf{lim_2}$ and $\sf{LPO^{\ast}}$ are incomparable, where $\sf{lim_2}:\subseteq 2^{\N}\ra2$ is the usual limit operation on $2^{\N}$. Therefore, we have the following corollary.
\begin{cor}
    $\sf{OG\mapsto\alpha}\not\geq_{\sf{W}}\sf{lim_2}$.
\end{cor}
Also, we can show the following.
\begin{prop}
    $\sf{OG\mapsto\alpha\epsilon}\not\geq_{\sf{W}}\sf{lim_2}\times\sf{lim_2}$.
\end{prop}

\begin{proof}
    It suffices to show that $\mathsf{WF}\times\mathsf{Min}\not\geq_{\sf{W}}\sf{lim_2}\times\sf{lim_2}$. Assume otherwise. We first fix forward and backward functionals $\Phi$, $\Psi$ witnessing $\mathsf{WF}\times\mathsf{Min}\geq_{\sf{W}}\sf{lim_2}\times\sf{lim_2}$. Let $(p_0,q_0)$ be a pair of infinite strings, where $p_0,q_0$ only consist of zeros. Let $\Phi(p_0,q_0)=(T_0,a_0)$ and $(\epsilon_0,b_0)$ be the output of $(\mathsf{WF}\times\mathsf{Min})(T_0,a_0)$. Then, $\Psi((p_0,q_0)\oplus (\epsilon_0,b_0))=(0,0)$. Let $u_0$ be the use of $\Psi$ with oracle $(p_0,q_0)\oplus(\epsilon_0,b_0)$. Let $v_0$ be the use of $\Phi$ for $b_0$ to appear for the first time in $a_0$. Let $w_0=\max\lb u_0,v_0\rb$. Note that $w_0$ is defined in this way to ensure that given any extensions $p'$, $q'$ of $p_0\upharpoonright w_0$, $q_0\upharpoonright w_0$, the minimum digit of $a'$ is at most $b_0$ where $(T',a')=\Phi(p',q')$ since $b_0$ has already appeared in $a'$. The idea is to continue locking in smaller and smaller upper bounds of the minimum digits, while forcing the outputs of the backward functional $\Psi$ to change.
    \par
    Let $q^0_1=q_0$ and $p_1^0$ be the infinite string starting with $p_0\upharpoonright w_0$ followed by all ones, and $(\epsilon_1^0,b_1^0)$ be the output of $\mathsf{WF}\times\mathsf{Min}$ with input $(T_1^0,a_1^0)=\Phi(p_1^0,q_1^0)$. Then, $\Psi((p_1^0,q_1^0)\oplus (\epsilon_1^0,b_1^0))=(1,0)$. Let $u_1^0$ be the use of $\Psi$ with oracle $(p_1^0,q_1^0)\oplus(\epsilon_1^0,b_1^0)$. Let $w_1^0=\max\lb w_0,u_1^0\rb$. Let $p_1^1$ be the infinite string starting with $p_1^0\upharpoonright w_1^0$ followed by all zeros, and $q_1^1$ be the infinite string starting with $q_1^0\upharpoonright w_1^0$ followed by all ones. Let $(\epsilon_1^1,b_1^1)$ be the output of $\mathsf{WF}\times\mathsf{Min}$ with input $(T_1^1,a_1^1)=\Phi(p_1^1,q_1^1)$. Then, $\Psi((p_1^1,q_1^1)\oplus (\epsilon_1^1,b_1^1))=(0,1)$.
    \par    
    Notice that $b_1^0,b_1^1\leq b_0$ since $b_0$ has appeared in $a_1^0$ and $a_1^1$. Also, there are only two possible values for $\epsilon_0,\epsilon_1^0,\epsilon_1^1$, i.e., $0$ and $1$. Then, $b_1^0$, $b_1^1$, and $b_0$ have at least two different possible values (so that $(\epsilon_1^0,b_1^0)$, $(\epsilon_1^1,b_1^1)$, and $(\epsilon_0,b_0)$ have more than two different values) since $\Psi$ has three different outputs. Therefore, at least one of $b_1^0$ and $b_1^1$ is strictly smaller than $b_0$. We let $b_1$ be one such number and define $p_1,q_1,a_1,T_1,\epsilon_1,u_1,v_1,w_1$ correspondingly.
    \par
    Then, we can define $q^0_2=q_1$ and $p_2^0$ to be the infinite string starting with $p_1\upharpoonright w_1$ followed by a tail of infinite $1-\mathsf{lim_2}(p_1)$'s, and $(\epsilon_2^0,b_2^0)$ be the output of $\mathsf{WF}\times\mathsf{Min}$ with input $(T_2^0,a_2^0)=\Phi(p_2^0,q_2^0)$. Let $u_2^0$ be the use of $\Psi$ with oracle $(p_2^0,q_2^0)\oplus(\epsilon_2^0,b_2^0)$. Let $w_2^0=\max\lb w_1,u_2^0\rb$. Let $p_2^1$ be the infinite string starting with $p_2^0\upharpoonright w_2^0$ followed by a tail of infinite $1-\mathsf{lim_2}(p_2^0)$'s, and $q_2^1$ be the infinite string starting with $q_2^0\upharpoonright w_2^0$ followed by a tail of infinite $1-\mathsf{lim_2}(p_2^0)$'s. Let $(\epsilon_2^1,b_2^1)$ be the output of $\mathsf{WF}\times\mathsf{Min}$ with input $(T_2^1,a_2^1)=\Phi(p_2^1,q_2^1)$. Similarly, we can define $b_2$ to be one of $b_2^0,b_2^1$ that is smaller than $b_1$.
    \par
    If we repeat this process, we will get an infinite strictly descending sequence of finite numbers. This is a contradiction. 
\end{proof}
\par
Notice that the proof can be modified slightly to show that $\mathsf{WF}\times\mathsf{Min}\not\geq_{\sf{W}}\sf{lim_3}$. Also, any problem with codomain $\lb0,1\rb$ can replace $\mathsf{WF}$ since this is the only property of $\mathsf{WF}$ used in the proof.
\begin{cor}\label{reverse}
    $\mathsf{C_{\N}}\not\leq\mathsf{OG\mapsto \alpha\epsilon}$.
\end{cor}
\begin{proof}
    Brattka, Gherardi, and Marcone showed that $\mathsf{C_{\N}}\equiv_{\sf{W}}\mathsf{lim_{\N}}>_{\sf{W}}\mathsf{lim_3}$ in \cite{lim2}. Also, notice that $\mathsf{C_{\N}}$ is a first-order problem. If $\mathsf{C_{\N}}\leq_{\sf{W}}\mathsf{OG\mapsto \alpha\epsilon}$, then $\mathsf{lim_3}<_{\sf{W}}\mathsf{C_{\N}}\leq_{\sf{W}}\mathsf{^1OG\mapsto \alpha\epsilon}\equiv_{\sf{W}}\mathsf{WF}\times\mathsf{Min}$. Contradiction.
\end{proof}
\par
Lastly, we would like to see whether it makes a difference if we restrict the order type of the group in the input to have the $\Q$-part or not have the $\Q$-part. We denote such problems $\mathsf{OG\mapsto\alpha1}$ and $\mathsf{OG\mapsto\alpha0}$, respectively.
\par
For the first-order part of $\mathsf{OG\mapsto\alpha0}$, notice that the proof of Proposition \ref{5} has shown that it is equivalent to $\mathsf{Min}$ as well. Therefore, we have the following proposition.

\begin{prop}
    $\sf{^1OG\mapsto\alpha0}\equiv_{\sf{W}}\sf{Min}$.
\end{prop}

\begin{cor}
    $\sf{OG\mapsto\alpha0}<_{\sf{W}}\sf{OG\mapsto\alpha1}\equiv_{\sf{W}}\sf{OG\mapsto\alpha}<_{\sf{W}}\sf{OG\mapsto\alpha\epsilon}$.
\end{cor}

\begin{proof}
    First, we show that $\sf{OG\mapsto\alpha1}\geq_{\sf{W}}\sf{OG\mapsto\alpha}$. Given a countable ordered group $G_0$, the forward functional $\Phi$ can build a group $G_1$ that is the direct product of $G_0$ and the additive group $\Q$ of rational numbers. When building $G_1$, $\Phi$ also orders the elements in $G_1$ such that for any $(g_0,q_0),(g_1,q_1)\in G_0\times\Q$, $(g_0,q_0)<_{G_1}(g_1,q_1)$ when $q_0<_{\Q}q_1$ or $q_0=q_1\wedge g_0<_{G_0}g_1$. Notice that for any $(g_0,q_0)<_{G_1}(g_1,q_1)$ and $(g_2,q_2)\in G_1$, $(g_0g_2,q_0+q_2)<_{G_1}(g_1g_2,q_1+q_2)$ and $(g_2g_0,q_2+q_0)<_{G_1}(g_2g_1,q_2+q_1)$. So, $G_1$ is an ordered group. Its order type is $\Z^{\alpha}\Q$ when $G_0$ has order type $\Z^{\alpha}$ or $\Z^{\alpha}\Q$ since $\Z^{\alpha}\Q\Q\cong\Z^{\alpha}\Q$. The reduction follows since we can let the backward functional output the copy of $\alpha$ in its input.
    \par
    Next, $\sf{OG\mapsto\alpha}\not\leq_{\sf{W}}\sf{OG\mapsto\alpha0}$ since $\alpha$ can be $\omega_1^{\mathrm{CK}}$ in the former given a computable input and a computable output of the latter cannot be used to compute a copy of it. Kihara, Marcone, and Pauly showed that $\mathsf{WF}\not\leq_{\sf{W}}\mathsf{C_{\N^{\N}}}$ in \cite{WFCNN}. It follows that $\mathsf{WF}\not\leq_{\sf{W}}\mathsf{LPO^{\ast}}$. By this fact, and Propositions \ref{W} and \ref{LW}, we conclude that $\sf{OG\mapsto\alpha\epsilon}\not\leq_{\sf{W}}\sf{OG\mapsto\alpha}$. Otherwise, $\sf{LPO^{\ast}}\times\sf{WF}\equiv_{\sf{W}}\sf{^1OG\mapsto\alpha\epsilon}\leq_{\sf{W}}\sf{OG\mapsto\alpha}$, contradicting $\sf{^1OG\mapsto\alpha}\equiv_{\sf{W}}\sf{LPO^{\ast}}$.
\end{proof}

\begin{prop}
    $\sf{OG\mapsto\alpha0}\leq_{\sf{W}}\sf{lim}$.
\end{prop}

This follows from Lemmas \ref{Arch} and \ref{G} that $\mathrm{Arch}(G)$ is a copy of $\alpha$ and uniformly computed from $G'$, and $G'$ is given by $\mathsf{lim}$. 
\par
The proposition above cannot be reversed by Corollary \ref{reverse}. However, $\sf{OG\mapsto\alpha0}$ does have some computational power. We introduce an intermediate problem that transforms a $\Sigma_3^0$-question into a $\Sigma_2^0$-question similar to the problem $\sf{\chi_{\Pi_2^0\ra\Pi_1^0}}$ in \cite{jump} defined by Andrews et al.

\begin{defn}
    We define $\sf{\chi_{\Sigma_3^0\ra\Sigma_2^0}}:\N^{\N}\rightrightarrows\N^{\N}$ as:
    \begin{align*}
        \mathsf{\chi_{\Sigma_3^0\ra\Sigma_2^0}}(p)=\lb q\in\N^{\N}:(\exists a\forall i>a\exists j)\lk p(\langle j,i\rangle)=1\rk\Leftrightarrow(\exists a\forall b>a)\lk q(b)=1\rk\rb.
    \end{align*}
\end{defn}

\begin{prop}\label{chi}
    $\sf{OG\mapsto\alpha0}\geq_{\sf{W}}\sf{\chi_{\Sigma_3^0\ra\Sigma_2^0}}$.
\end{prop}

\begin{proof}
    Let $p$ be an input of $\sf{\chi_{\Sigma_3^0\ra\Sigma_2^0}}$. The forward functional $\Phi$ builds an ordered group $G$ of order type $\Z^{\omega}$ or $\Z^n$ for some finite number $n$ depending on $p$. We say that for any $j$, a digit $p(\langle j,i\rangle)$ of $p$ is in the $i$th column.
    \par
    Before $\Phi$ reads any digit of $p$, it puts a generator in $G$ and identifies it with $(1)$ in order type $\Z$.
    \par
    Suppose at stage $s$, $\Phi$ reads a new digit $p(d_s)$ of $p$, and $\Phi$ has put elements into $G$ so that $G$'s elements are identified with elements in order type $\Z^k$, $k\geq 1$. At this stage $s$, $\Phi$ operates in two substages. In substage $0$, there are two cases.
    \par
    Case $1$: when $p(d_s)\not=1$ is the first digit in column $i$ read by $\Phi$, $\Phi$ puts a new largest generator into $G$ and identifies it with $(\underbrace{0,\ldots,0}_{k\text{ many}},1)$. Also, any $g\in G$ identified with $(a_0,\ldots,a_{k-1})$ is now identified with $(a_0,\ldots,a_{k-1},0)$. That is, all elements in $G$ are identified with elements in the order type $\Z^{k+1}$.
    \par
    Case $2$: when $p(d_s)=1$ is not the first digit in column $i$ read by $\Phi$, but the first $1$ read by $\Phi$ in this column, $\Phi$ reinterprets $G$'s elements. Suppose that the generator put into $G$ when $\Phi$ first saw a digit in column $i$ is now identified with $(\underbrace{0,\ldots,0}_{n\text{ many}},1,0,\ldots0)$, where $n\geq 1$. There exists a large enough $l$ similar to the proof of Proposition \ref{5} such that for any $g\in G$ identified with $(a_0,\ldots,a_{n-1},a_{n},a_{n+1},\ldots,a_{k-1})$, we can map it to $(a_0,\ldots,a_{n-1}+a_{n}\cdot l,a_{n+1},\ldots,a_{k-1})$ consistently with the current presentation. We now interpret each element of $G$ as an element in $\Z^{k-1}$. 
    \par
    In substage $1$, assume that $G$'s elements are identified with elements in order type $\Z^{k'}$, $k'\geq 1$. In this substage, 
    $\Phi$ puts new elements into $G$. Such elements are generated by existing generators and identified with elements in $\Z^{k'}$. Here, $\Phi$ puts enough elements into $G$ so that for any $(a_0,\ldots,a_{k'-1})$, where $|a_j|\leq s$ when $0\leq j\leq k'-1$, it is identified with a group element. Here, $k'\geq 1$ because a generator was put into $G$ before $\Phi$ reads any digit of $p$. This ensures that in later stages, the reinterpretation in substage $0$ can be done when $n=1$.
    \par
    After $\Phi$ reads $p(d_s)$ and acts accordingly, it moves on to read the next digit of $p$.
    \par
    In this way, when $p$ satisfies $(\exists a\forall i>a\exists j)\lk p(\langle j,i\rangle)=1\rk$, $\Phi$ produces a group of order type $\Z^{k+1}$ where $k$ is the number of $i$'s such that $(\forall j)\lk p(\langle j,i\rangle)\not= 1\rk$. Otherwise, $\Phi$ produces a group of order type $\Z^{\omega}$. 
    \par
    Given a copy $c\in\mathsf{OG\mapsto\alpha0}(\Phi(p))$, the backward functional $\Psi$ determines whether there is a largest element in $c$. If the $s$th bit in $c$ suggests that there is an element $c_s$ in $c$, $\Psi$ checks if all elements appeared before this bit are smaller than $c_s$. If so, let $q(s)=0$. If not, or the $s$ bit does not give an element in $c$, let $q(s)=1$. Then, $q$ satisfies $(\exists a\forall b>a)\lk q(b)=1\rk$ if and only if $c$ is a copy of a finite number. This is because there are infinitely many $0$'s if and only if there is not a largest element in $c$.
\end{proof}

\begin{prop}\label{chi'}
    $\sf{C_{\N}}\not\geq_{\sf{W}}\sf{\chi_{\Sigma_3^0\ra\Sigma_2^0}}$.
\end{prop}

\begin{proof}
    Assume otherwise. Let $\Phi,\Psi$ be the forward and backward functionals, respectively. Let $\mathrm{Cof}=\lb e:\mathrm{dom}(\phi_e)\mathrm{\ is\ cofinite}\rb$. Given each $e\in\N$, we can build a $p_e\in\N^{\N}$ in the following way: $p_e(\langle j,i\rangle)=1$ when $\phi_{e,j}(i)\downarrow$ and $p_e(\langle j,i\rangle)=0$ otherwise. Then, $p_e$ satisfies $(\exists a\forall i>a\exists j)\lk p_e(\langle j,i\rangle)=1\rk$ if and only if $e\in\mathrm{Cof}$. Let $W_{i_e}$ be the c.e.\ subset such that $\mathsf{C_{\N}}(\Phi(p_e))=\N\setminus W_{i_e}$. Notice that a $0'$ oracle can be used to produce an $n\not\in W_{i_e}$. Also, $0''$ can determine whether $q_e=\Psi(p_e,n)$ satisfies $(\exists a\forall b>a)\lk q_e(b)=1\rk$ or not since $0''$ is $\Sigma_2^0$-complete. Therefore, $0''$ can determine whether each $e\in\mathrm{Cof}$ or not. This contradicts the fact that $\mathrm{Cof}$ is $\Sigma_3^0$-complete.
\end{proof}
Propositions \ref{chi} and \ref{chi'} give us the following corollary.
\begin{cor}\label{last}
    $\sf{C_{\N}}\not\geq_{\sf{W}}\sf{OG\mapsto\alpha0}$.
\end{cor}

Finally, we show that $\mathsf{OG\mapsto \alpha1}$, $\mathsf{OG\mapsto \alpha}$, and $\mathsf{OG\mapsto \alpha0}$ are incomparable to $\mathsf{WF}$.

\begin{prop}
    $\sf{WF}\mathbin{|_{\sf{W}}}\sf{OG\mapsto\alpha0}$, $\sf{WF}\mathbin{|_{\sf{W}}}\sf{OG\mapsto\alpha1}$, and $\sf{WF}\mathbin{|_{\sf{W}}}\sf{OG\mapsto\alpha}$. 
\end{prop}

\begin{proof}
    It suffices to show that $\sf{WF}\mathbin{|_{\sf{W}}}\sf{LPO^{\ast}}$ because $\sf{WF}$ is a first-order problem and $\sf{LPO^{\ast}}$ is the first-order part of $\sf{OG\mapsto\alpha0}$, $\sf{OG\mapsto\alpha1}$, and $\sf{OG\mapsto\alpha}$. 
    \par
    Given an infinite binary string that has a limit, we can build a tree that extends to a single child on each layer if and only if we find the next $0$ in the string. The tree is ill-founded if and only if the limit is $0$. So, $\mathsf{lim_2}\leq_{\sf{W}} \mathsf{WF}$. Then, by the incomparability of $\mathsf{lim_2}$ and $\mathsf{LPO^{\ast}}$, $\mathsf{WF}\not\leq_{\sf{W}}\sf{LPO^{\ast}}$. 
    \par
    For the other direction, it suffices to show $\mathsf{LPO\times LPO}\not\leq_{\mathsf{W}}\mathsf{WF}$. Assume otherwise. We fix the forward and backward functionals $\Phi$, $\Psi$. Let $(p_0,q_0)$ be a pair of infinite strings, where $p_0, q_0$ only consist of ones. Let $\Phi(p_0,q_0)=T_0$ and $\epsilon_0$ be the output of $\mathsf{WF}(T_0)$. Then, $\Psi((p_0,q_0)\oplus \epsilon_0)=(1,1)$. Let $u_0$ be the use of $\Psi$ with oracle $(p_0,q_0)\oplus \epsilon_0$. Let $q_1=q_0$ and $p_1$ be the infinite string starting with $p_0\upharpoonright u_0$ followed by all zeros, and $\epsilon_1$ be the output of $\mathsf{WF}$ with input $T_1=\Phi(p_1,q_1)$. Then, $\Psi((p_1,q_1)\oplus \epsilon_1)=(0,1)$. Let $u_1'$ be the use of $\Psi$ with oracle $(p_1,q_1)\oplus \epsilon_1$. Let $u_1=\max\lb u_0,u_1'\rb$. Let $p_2=p_1$ and $q_2$ be the infinite string starting with $q_1\upharpoonright u_1$ followed by all zeros. Let $\epsilon_2$ be the output of $\mathsf{WF}$ with input $T_2=\Phi((p_2,q_2))$. Then, $\Psi((p_2,q_2)\oplus \epsilon_2)=(0,0)$. Notice that $\epsilon_0$ and $\epsilon_1$ are different since $(1,1)=\Psi((p_0,q_0)\oplus \epsilon_0\upharpoonright u_0)=\Psi((p_1,q_1)\oplus \epsilon_1\upharpoonright u_0)=(0,1)$ otherwise. Similarly, we can show that $\epsilon_0,\epsilon_1,\epsilon_2$ are pairwise distinct. However, there are only two possible values for $\epsilon_0,\epsilon_1,\epsilon_2$, i.e., $0$ and $1$. Contradiction.
\end{proof}

\bibliographystyle{plainurl}
\bibliography{references} 

\end{document}